\documentclass[a4paper,11pt]{article}
\usepackage{amsmath}
\usepackage{amssymb}
\usepackage{graphicx}
\usepackage{parskip}
\usepackage{amsthm}
\usepackage{fancyhdr}
\usepackage{dsfont}
\usepackage{algorithm}
\usepackage{algorithmic}
\usepackage{hyperref}
\usepackage{bm}
\usepackage[authoryear]{natbib}

\newtheorem{thm}{Theorem}
\newtheorem{prop}{Proposition}
\newtheorem{lem}{Lemma}
\newtheorem{cor}{Corollary}
\theoremstyle{definition}
\newtheorem{defn}{Definition}

\newtheorem{rmk}{Remark}

\newcommand{\M}{ \mathcal{ M } }

\newcommand{\N}{ \mathbb{ N } }
\newcommand{\E}{ \mathbb{ E } }
\newcommand{\R}{ \mathbb{ R } }
\newcommand{\bbZ}{ \mathbb{ Z } }
\newcommand{\Prb}{ \mathbb{ P } }

\newcommand{\x}{ \mathbf{ x } }
\newcommand{\y}{ \mathbf{ y } }
\newcommand{\Y}{ \mathbf{ Y } }
\newcommand{\bfk}{ \mathbf{ k } }
\newcommand{\X}{ \mathbf{ X } }

\newcommand{\z}{ \mathbf{ z } }
\newcommand{\Z}{ \mathbf{ Z } }
\newcommand{\W}{ \mathbf{ W } }
\newcommand{\KL}{ \operatorname{KL} }

\newcommand{\bxi}{ \bm{ \xi } }
\newcommand{\bfzeta}{ \bm{ \zeta } }
\newcommand{\bflambda}{ \bm{ \lambda } }

\title{Consistency of Bayesian nonparametric inference for discretely observed jump diffusions}

\author{Jere Koskela \\
	\texttt{j.koskela@warwick.ac.uk}\\
	\small Department of Statistics \\
	\small University of Warwick \\
	\small Coventry CV4 7AL \\
	\small United Kingdom
	\and
	Dario Span\`{o} \\
	\texttt{d.spano@warwick.ac.uk}\\
	\small Department of Statistics \\
	\small University of Warwick \\
	\small Coventry CV4 7AL \\
	\small United Kingdom
	\and
	Paul A. Jenkins \\
	\texttt{p.jenkins@warwick.ac.uk}\\
	\small Departments of Statistics and Computer Science\\
	\small University of Warwick \\
	\small Coventry CV4 7AL \\
	\small United Kingdom
}
\date{\today}

\begin{document}

\maketitle

\begin{abstract}
We introduce verifiable criteria for weak posterior consistency of Bayesian nonparametric inference for jump diffusions with unit diffusion coefficient and uniformly Lipschitz drift and jump coefficients in arbitrary dimension.
The criteria are expressed in terms of coefficients of the SDEs describing the process, and do not depend on intractable quantities such as transition densities.
We also show that priors built from discrete nets, wavelet expansions, and Dirichlet mixture models satisfy our conditions.
This generalises known results by incorporating jumps into previous work on unit diffusions with uniformly Lipschitz drift coefficients.
\end{abstract}

\section{Introduction}\label{introduction}

Jump diffusions are a broad wide class of stochastic processes encompassing systems undergoing deterministic mean-field dynamics, microscopic diffusion and macroscopic jumps.
In this paper we let $\X := ( \X_t )_{ t \geq 0 }$ denote a unit jump diffusion, which can be described as a solution to a stochastic differential equation of the form
\begin{align}
d\X_t = b( \X_t ) dt + d\W_t + c( \X_{ t- }, d\Z_t ) \label{sde}
\end{align}
on a domain $\Omega \subseteq \R^d$ given an initial condition $\X_0 = \x_0$, coefficients $b : \Omega \mapsto \R^d$ and $c : \Omega \times \R_0^d \mapsto \R_0^d$, a $d$-dimensional Brownian motion  $( \W_t )_{ t \geq 0 }$ and a pure jump L\'evy process $( \Z_t )_{ t \geq 0 }$ on $\R_0^d := \R^d \setminus \{ \mathbf{ 0 } \}$ with L\'evy measure $M( d\z )$ satisfying
\begin{equation*}
\int_{ \R_0^d } ( \| \z \|_2^2 \wedge 1 ) M( d\z ) < \infty
\end{equation*}
The notation $\| \cdot \|_{ p, \rho }$ denotes the $L^p( \rho )$-norm, where the Lebesgue measure is meant whenever the measure $\rho$ is omitted.

Jump diffusions are used as models across a broad spectrum of applications, such as economics and finance \citep{Merton76, Aase87, Bardhan93, Chen05, Filipovic07}, biology \citep{Kallianpur92, Kallianpur94, Bertoin03, Birkner09a} and engineering \citep{Au82, Bodo87}.
They also contain many important families of stochastic processes as special cases, including diffusions and L\'evy processes.
\vskip 11pt
\begin{rmk}
In the exposition above, the processes $\X$, $\W$ and $\Z$ all share a common dimension.
This restriction is not necessary for any of the results in the paper, and has been introduced purely for readability of notation.
\end{rmk}

Under regularity conditions summarised in the next section, jump diffusions are recurrent, ergodic Feller-Markov processes with transition densities $p_t( \x, \y )d\y$ and a unique stationary density $\pi( \x ) d\x$ with respect to the $d$-dimensional Lebesgue measure.
Under such conditions the procedure of Bayesian inference can be applied to infer the coefficients of the jump diffusion based on observations taken at discrete times.
In this paper we focus on joint inference of the drift function $b$ and the family of L\'evy measures $\nu( \x, d\z ) := M( c^*( \x, d\z ) )$, where $c^*( \x, \cdot )$ denotes the pull-back of $c( \x, \cdot)$:
\begin{equation*}
c^*( \x, d\z ) := \{ \y \in \R_0^d : c( \x, \y ) \in d\z \}.
\end{equation*}
We abuse terminology and refer to the collection of measures $\nu( \x, \cdot)$ as a L\'evy measure for the remainder of the paper.
Inference of the L\'evy measure will refer to inference of $\nu$, assuming that neither $c$ nor $M$ is known.

More precisely, let $\Theta$ denote a set of pairs $( b,\nu )$, and let $\Pi$ denote a prior distribution on $( \Theta, \mathcal{ B }( \Theta ) )$, where $\mathcal{ B }( \Theta )$ is the Borel $\sigma$-algebra.
Let $\x_{ 0 : n } = ( \x_0, \x_{ \delta }, \ldots, \x_{ \delta n } )$ denote a time series of observations sampled from a stationary jump diffusion $\X$ at fixed separation $\delta$.
The object of interest is the posterior distribution, which can be expressed as
\begin{equation*}
\Pi( A | \x_{ 0 : n } ) := \frac{ \int_A \pi^{ b, \nu }( \x_0 ) \prod_{ i = 1 }^n p_{ \delta }^{ b, \nu }( \x_{ i - 1}, \x_i ) \Pi( db, d\nu ) }{ \int_{ \Theta } \pi^{ b, \nu }( \x_0 ) \prod_{ i = 1 }^n p_{ \delta }^{ b, \nu }( \x_{ i - 1}, \x_i ) \Pi( db, d\nu ) }
\end{equation*}
for measurable sets $A \in \mathcal{ B }( \Theta )$.
In the Bayesian setting, the posterior encodes all the available information for inferential purposes.
The restriction to unit diffusion coefficients implicit in \eqref{sde} is a strong assumption in dimension $d > 1$, though some models which fail to satisfy it outright can still be treated via the Lamperti transform \citep{Ait-Sahalia08}.
We will outline this procedure briefly in Section \ref{jump_diff}.

A typical approach to practical Bayesian inference is to choose $\Theta$ comprised of parametric families of drift functions and L\'evy measures, and fit these parameters to data.
However, the natural parameter spaces for jump diffusions are spaces of functions and measures, which are infinite dimensional and cannot be represented in terms of finitely many parameters without significant loss of modelling freedom.
Nonparametric Bayesian inference can be thought of as inference of infinitely many parameters, and retains much of the modelling freedom inherent in the class of jump diffusions.

A natural and central question is whether the Bayes procedure is \emph{consistent}, that is, whether the posterior concentrates on a neighbourhood of the parameter space which specifies the ``true" dynamics generating the data as the number of observations grows.
If $( b_0, \nu_0 ) \in \Theta$ denotes the data generating drift and L\'evy measure, consistency can be expressed as $\Pi( U_{ b_0,\nu_0 }^c | \x_{ 0 : n } ) \rightarrow 0$ as $n \rightarrow \infty$, where $U_{ b_0,\nu_0 }$ is any open neighbourhood of $( b_0,\nu_0 )$.
In the context of \emph{non-identifiable} problems, where the prior supports multiple collections of parameters which the data cannot distinguish from the ``true" dynamics, consistency is instead defined as $\Pi( ( U_{ b_0, \nu_0 }^{ TI } )^c | \x_{ 0 : n } ) \rightarrow 0$ as $n \rightarrow \infty$, where $U_{ b_0, \nu_0 }^{ TI }$ is the set of points that are topologically indistinguishable from $( b_0, \nu_0 )$ in the parameter space.

Whether or not Bayesian posterior consistency holds in the nonparametric setting is an intricate question, and depends on subtle ways on the prior $\Pi$ and the topology endowed on $\Theta$ \citep{Diaconis86}.
A further difficulty in the present context is the fact that stationary and transition densities of jump diffusions are intractable in practically all cases of interest, so that usual conditions for posterior consistency are difficult to verify.
These difficulties were recently overcome for discretely observed, one-dimensional unit diffusions under restrictive conditions on the drift function \citep{vanderMeulen13}, and a multidimensional generalisation was presented in \citep{Gugushvili14}.
Both results rely on martingale arguments developed by Ghosal and Tang for Markov processes with tractable transition probabilities \citep{Ghosal06, Tang07}.
A Bayesian analysis of continuously observed one dimensional diffusions has also been conducted under various setups \citep{vanderMeulen06, Panzar09, Pokern13}, and a review of Bayesian methods for one dimensional diffusions is provided by \citep{vanZanten13}.
Similar developments have also been made for frequentist drift estimation from discrete observations, both for one dimensional unit diffusions \citep{Jacod00, Gobet04, Comte07} and their multi-dimensional generalisations \citep{Dalalyan07, Schmisser13}.

The main result of this paper is consistency of Bayesian nonparametric joint inference of drift functions and L\'evy measures in arbitrary dimension under verifiable conditions on the prior.
This generalises the result of \citep{Gugushvili14} by incorporating discontinuous processes with jumps.
We also show that our consistency conditions are satisfied by a prior $\Pi$ constructed by specifying independent discrete net or wavelet expansion priors for $b$ and $c$, and a further independent Dirichlet mixture model prior for $M$.
The key results enabling this generalisation are a generalised Girsanov-type change of measure theorem for jump diffusions \citep{Cheridito05}, and a coupling method for establishing regularity of semigroups \citep{Wang10}.

The rest of the paper is organised as follows.
In Section \ref{jump_diff} we introduce the jump diffusion processes in finite dimensional domains and necessary regularity conditions.
In Section \ref{consistency} we define the inference problem under study, and state and prove the corresponding consistency result.
In Section \ref{prior} we introduce the discrete net prior, and show that it satisfies our consistency conditions.
Section \ref{discussion} concludes with a discussion.

\section{Jump diffusions}\label{jump_diff}

A general time-homogeneous, $d$-dimensional jump diffusion $\Y := ( \Y_t )_{ t \geq 0 }$ is the solution of a stochastic differential equation of the form
\begin{equation*}
d\Y_t = b( \Y_t ) dt + \sigma( \Y_t )d\W_t + c( \Y_{ t- }, d\Z_t ),
\end{equation*}
where $\sigma : \Omega \mapsto \R^{ d \times d }$ and the other coefficients are as in \eqref{sde}.
The implicit assumption in \eqref{sde} of $\sigma \equiv 1$ is restrictive in dimensions $d > 1$.
Processes which do not have unit diffusion coefficient can be dealt with provided they lie in the domain of the Lamperti transform \citep{Ait-Sahalia08}, i.e.~if there exists a mapping  $q : \Y \mapsto \X$ such that $\X$ is of the form \eqref{sde}.
Such transforms exist for any non-degenerate process in one dimension, but only rarely in higher dimensions.
Sufficient conditions for the Lamperti transform to be well defined are non-singularity of $\sigma$ and the following symmetry condition \citep{Yu07, Ait-Sahalia08}:
\begin{equation}\label{Lamperti_condition}
\frac{ \partial ( \sigma^{ -1 } )_{ i j }( \x ) }{ \partial x_k }  =\frac{ \partial ( \sigma^{ -1 } )_{ i k }( \x ) }{ \partial x_j } \text{ for all } i, j, k \in \{ 1, \ldots, d \}.
\end{equation}
We note also that the Lamperti transform cannot be constructed from discrete data, so that in any case $\sigma$ must be known a priori.
While restrictive, this assumption cannot be relaxed without fundamental changes to the method of proof of consistency and already arises in the simpler case of diffusions without jumps \citep{vanderMeulen13, Gugushvili14}.

The following proposition summarises the necessary regularity assumptions for existence and uniqueness of Feller-Markov jump diffusions with transition densities and a unique stationary density:
\vskip 11pt
\begin{prop}\label{existence}
Assume that $c( \cdot, 0 ) \equiv 0$, and that there exist constants $C_1, C_2, C_3, C_4 > 0$ such that
\begin{align}
\| b( \mathbf{ x } ) - b( \mathbf{ y } ) \|_2^2 + \int_{ \R_0^d } \| c( \mathbf{ x }, \z ) - c( \mathbf{ y }, \z ) \|_2^2 M( d\z ) &\leq C_1 \| \mathbf{ x } - \mathbf{ y } \|_2^2 \label{lip} \\
\| c( \mathbf{ x }, \z ) - c( \mathbf{ x }, \bxi ) \|_2^2 &\leq C_2 \| \z - \bxi \|_2^2 \label{lip2}\\
\text{For every } \mathbf{ x } \in \Omega : \| \mathbf{ x } \|_2 > C_3 \text{ the following holds: } \mathbf{ x } \cdot b( \mathbf{ x } ) &\leq - C_4 \| \mathbf{ x } \|_2^2 \label{drift} \\
\int_{ \mathbb{ R }_0^d : \| \z \|_2 > 1 } \| \mathbf{ z } \|_2^2 M( dz ) &< \infty. \label{m_moment}
\end{align}
Then \eqref{sde} has a unique, ergodic weak solution $\X$ with the Feller and Markov properties.
Furthermore, $\X$ has a unique stationary density $\pi^{ b, \nu }( \mathbf{ x } ) d\mathbf{ x }$ with a finite second moment, and the associated semigroup $P_t^{ b, \nu }$ has transition densities $p_t^{ b, \nu }( \mathbf{ x }, \mathbf{ y } ) d\mathbf{ y }$.
\end{prop}
\begin{proof}
Existence and uniqueness of $\X$ are obtained from \eqref{lip}, as well as the linear growth bounds implied by Lipschitz continuity, by Theorem 6.2.9 of \citep{Applebaum04}.
Theorem 6.4.6 of \citep{Applebaum04} gives the Markov property under the same conditions.
Finally, the corollary in Appendix 1 of \citep{Kolokoltsov04} yields the Feller property.
In turn, the Feller property and the fact that $\log( 1 + \| \bxi \|_2 )^{ -1 } \| \bxi \|_2^2 \rightarrow \infty$ as $\| \bxi \|_2 \rightarrow \infty$ mean that the hypotheses of Theorem 1.1 of \citep{Schilling13} are fulfilled, so that $\X$ has bounded transition densities with respect to the Lebesgue measure.

Existence and uniqueness of $\pi^{ b, \nu }$, as well as ergodicity of $\X$ will follow from Theorem 2.1 of \citep{Masuda07}, the hypotheses of which will now be verified.
Along with $c( \cdot, 0 ) \equiv 0$, conditions \eqref{lip} and \eqref{lip2} above imply  Assumption 1 of \citep{Masuda07}.
Now, for every $u \in ( 0, 1 )$ let 
\begin{equation*}
b^u( \mathbf{ x } ) := b( \mathbf{ x } ) - \int_{ u < \| \z \|_1 \leq 1 } c( \mathbf{ x }, \z ) M( d\z ).
\end{equation*}
Assumption 2(a)' of \citep{Masuda09} requires $\X$ to admit bounded transition densities, and the diffusion which solves
\begin{equation*}
d\X_t^u = b^u( \X_t^u ) dt + \sigma( \X_t^u ) d\W_t
\end{equation*}
to be irreducible for each $u > 0$.
Boundedness of the transition density of $\mathbf{ X }$ was established above, and irreducibility of $\mathbf{ X }^u$ holds because $\sigma \equiv 1$ by Theorem 2.3 of \citep{Stramer97}.

Next we verify Assumptions 3 and 3* of \citep{Masuda07} by checking the conditions of Lemma 2.4' of \citep{Masuda09}.
The diffusion coefficient is constant, and hence $o( \| x \|_2^{ 1 - q / 2 } )$ for any $q \in ( 0, 2 )$.
Condition \eqref{m_moment} is the corresponding hypothesis of \citep{Masuda09}, and both $\| \mathbf{ x } \|_2^{ q - 2 } \mathbf{ x } \cdot b( \mathbf{ x } ) \rightarrow - \infty$ and $\| \mathbf{ x } \|_2^{ - 2 } \mathbf{ x } \cdot b( \mathbf{ x } ) \leq - C_4$ follow from \eqref{drift}.
Hence, Assumptions 3 and 3* of \citep{Masuda07} hold.
This yields ergodicity (and mixing) by Theorem 2.1 of \citep{Masuda07}, and second moments of the stationary distribution (and exponential mixing) by Theorem 2.2 of \citep{Masuda07}.

It remains to show the invariant measure has a density.
By combining Proposition 5.1.9 and Theorem 5.1.8 of \citep{Fornaro04} it can be seen that invariant measures of irreducible strong Feller processes are equivalent to the associated transition probabilities, which is sufficient in this case.
Assumption 1 of \citep{Masuda07} and Assumption 2(a)' of \citep{Masuda09} imply irreducibility of $\X$ (c.f.~Claim 1 on page 42 of \citep{Masuda07}).
Condition \eqref{lip} guarantees the strong Feller property by Theorem 2.3 of \citep{Wang10}.
Hence the invariant measure has a density with respect to the transition densities, and thus also the Lebesgue measure.
This concludes the proof.
\end{proof}
\begin{rmk}
Assumption \eqref{lip} is central to the proof of our main result.
In contrast, assumptions \eqref{lip2}, \eqref{drift} and \eqref{m_moment} are only needed to ensure the conclusions of Proposition \ref{existence}.
\end{rmk}

We denote the law of $\X$ with drift function $b$, L\'evy measure $\nu$ and initial condition $\X_0 = \x$ by $\Prb_ { \x }^{ b, \nu }$ and the corresponding expectation by $\E_{ \x }^{ b, \nu }$.
Dependence on initial conditions is omitted when the stationary process is meant.

\section{Consistency for discrete observations}\label{consistency}

We begin by defining the topology and weak posterior consistency following the set up of \citep{vanderMeulen13}.
In addition to topological details, posterior consistency is highly sensitive to the support of the prior, which should not exclude the truth.
This is guaranteed by insisting that the prior places positive mass on all neighbourhoods of the truth, typically measured in terms of Kullback-Leibler divergence.
In our setting such a support condition is provided by \eqref{kl_cond} below.

We begin by setting out the necessary assumptions on the parameter space $\Theta$.
\vskip 11pt
\begin{defn}\label{prior_supp}
Let $\Theta \subseteq \{ ( b, \nu ) : b : \Omega \mapsto \R^d, \nu : \Omega \times \R_0^d \mapsto \R_+ \}$ denote a set of pairs of drift functions $b( \x )$ and L\'evy measures $\nu( \x, d\z ) := M( c^*( \x, d\z ) )$.
Suppose each pair satisfies the hypotheses of Proposition \ref{existence}, and that there exists a constant $C_5 > 0$ such that
\begin{equation}\label{m_mass}
M( \{ \z \in \R_0^d : \| \z \|_2 > 1 \} ) < C_5
\end{equation}
holds uniformly in $\Theta$.
In addition, let $\Theta$ be such that $\nu( \x, \cdot ) \sim \nu'( \x, \cdot )$, that
\begin{equation}\label{rn_bounds}
0 < \inf_{ \substack{ \x \in \Omega \\ \z \in \R_0^d } } \left\{ \frac{ d\nu' }{ d\nu }( \x, \z ) \right\} \leq \sup_{ \substack{ \x \in \Omega \\ \z \in \R_0^d } } \left\{ \frac{ d\nu' }{ d\nu }( \x, \z ) \right\} < \infty,
\end{equation}
and that either
\begin{enumerate}
\item $\nu( \x, \cdot )$ and $\nu'( \x, \cdot )$ are finite measures, or
\item there exists an open set $A$ containing the origin such that $\nu( \x, \cdot )|_A = \nu'( \x, \cdot )|_A$,
\end{enumerate}
for any pair $( \cdot, \nu), (\cdot, \nu' ) \in \Theta$ and for each $\x \in \Omega$.
\end{defn}
\vskip 11pt
\begin{rmk}
In effect, the conditions of Definition \ref{prior_supp} mean that the unit diffusion coefficient and the infinite intensity component of the L\'evy measure are known confounders of the joint inference problem for the drift function and the compound Poisson component of the L\'evy measure driving macroscopic jumps.
In particular, one of conditions 1.~or 2.~is needed to ensure finiteness of the second and third terms in \eqref{kl_cond}.
\end{rmk}

The topology under consideration is defined as in \citep{vanderMeulen13, Gugushvili14} by specifying a subbase determined by the semigroups $P_t^{ b, \nu }$.
For details about the notion of a subbase, and other topological concepts, see e.g.~\citep{Dudley02}.
Let $\M_f( \Omega )$ denote the set of finite measures supported on a set $\Omega$, and let $C_b( \Omega )$ denote the Banach space of continuous, bounded functions on $\Omega$.
\vskip 11pt
\begin{defn}\label{topology}
Fix a sampling interval $\delta > 0$ and a finite measure $\rho \in \M_f( \Omega )$ with positive mass in all non-empty, open sets.
For any $( b, \nu ) \in \Theta$, $\varepsilon > 0$ and $f \in C_b( \Omega)$ define  the set
\begin{equation*}
U_{ f, \varepsilon }^{ b, \nu } := \{ ( b', \nu' ) \in \Theta : \| P_{ \delta }^{ b', \nu' } f - P_{ \delta }^{ b, \nu } f  \|_{ 1, \rho } < \varepsilon \}.
\end{equation*}
A weak topology on $\Theta$ is generated by requiring that the family $\{ U_{ f, \varepsilon }^{ b, \nu } : f \in C_b( \Omega ), \varepsilon > 0, ( b, \nu ) \in \Theta \}$ is a subbase of the topology.
\end{defn}
The topology introduced in Definition \ref{topology} coincides with that used by both \citet{vanderMeulen13} for scalar diffusions, as well as \citet{Gugushvili14} for multi-dimensional diffusions.
In both of these works it was shown to separate points, so that the corresponding inference problems are identifiable.
The key tool was a one-to-one correspondence between drift functions and stationary laws of the corresponding diffusion process.
Because of jumps, such a correspondence does not hold in our setting in general, and thus we have the following, weaker separation result.
\vskip 11pt
\begin{lem}\label{ident_lemma}
The topology defined in Definition \ref{topology} is \emph{pre-regular}, that is for any pair $( b, \nu ) \neq ( b' \nu' ) \in \Theta$ with $P_{ \delta }^{ b', \nu' } \neq P_{ \delta }^{ b, \nu }$ there exists a function $f \in C_b( \Omega )$ and $\varepsilon > 0$ such that the neighbourhoods $U_{ f, \varepsilon }^{ b, \nu }$ and $U_{ f, \varepsilon }^{ b', \nu' }$ are disjoint.
\end{lem}
\begin{proof}
By assumption there exists a function $f \in C_b( \Omega )$ and a point $\x \in \Omega$ such that $P_{ \delta }^{ b, \nu } f( \x ) \neq P^{ b', \nu' } f( \x )$.
Thus, continuity implies that there exists a non-empty, open set $J$ on which the functions $P^{ b, \nu } f$ and $P^{ b', \nu' } f$ differ.
Since $\rho$ assigns positive mass to this set by assumption, we have $\| P^{ b, \nu } f - P^{ b', \nu' } f \|_{ 1, \rho } > \varepsilon > 0$ for some $\varepsilon > 0$.
\end{proof}
The content of Lemma \ref{ident_lemma} is that our topology can separate distinct pairs of drift and jump coefficients whenever they give rise to distinct transition semigroups across time interval $\delta > 0$, the time-spacing of the observations.
Obviously, no statistical procedure based on discrete data with time separation $\delta > 0$ can separate two pairs of coefficients for which these transition semigroups coincide, and hence this topology is well suited to our purposes.

We are now in a position to formally define posterior consistency, and state the main result of the paper.
\vskip 11pt
\begin{defn}\label{def_cons}
Let $\x_{ 0 : n } := ( \x_0, \ldots, \x_n )$ denote $n + 1$ samples observed at sampling times $0 , \delta, \ldots, \delta n$  from $\X$ at stationarity, i.e.~with initial distribution $\X_0 \sim \pi^{ b_0, \nu_0 }$.
Weak posterior consistency holds if $\Pi( ( U_{ b_0, \nu_0 }^{ TI } )^c | \x_{ 0 : n } ) \rightarrow 0$ with $\Prb^{ b_0, \nu_0 }$-probability 1 as $n \rightarrow \infty$, where $U_{ b_0, \nu_0 }^{ TI }$ is any open neighbourhood of the set of points that are topologically indistinguishable from $( b_0, \nu_0 ) \in \Theta$.
\end{defn}
\vskip 11pt
\begin{thm}\label{ts_consistency}
Let $\x_{ 0 : n }$ be as in Definition \ref{def_cons}, and suppose that the prior $\Pi$ is supported on a set $\Theta$ which satisfies the conditions of Definition \ref{prior_supp}, with the constant $C_1$ in \eqref{lip} holding uniformly in $\Theta$.
If 
\begin{align}
\Pi\Bigg( ( b, \nu ) \in \Theta &: \frac{ 1 }{ 2 } \E^{ b_0, \nu_0 }\left[ \| b_0( \X_0 ) - b( \X_0 ) \|_2^2 \right] \nonumber \\
&+ \frac{ 1 }{ 2 } \E^{ b_0, \nu_0 }\left[ \Big\| \int_{ B_0( 1 ) \setminus \{ 0 \} } \left( \frac{ d\nu_0 }{ d\nu }( \X_{ 0- }, \mathbf{ z } ) - 1 \right) \z \nu( \X_{ 0- }, d\mathbf{ z } ) \Big\|_2^2 \right]  \nonumber \\
&+ \E^{ b_0, \nu_0 }\left[ \int_{ \R_0^d } \left| \log\left( \frac{ d\nu_0 }{ d\nu }( \X_{ 0- }, \mathbf{ z } ) \right) - \frac{ d \nu_0 }{ d \nu }( \X_{ 0- }, \mathbf{ z } ) + 1 \right| \nu_0( \X_{ 0- }, d\mathbf{ z } ) \right] < \varepsilon \Bigg) > 0 \label{kl_cond}
\end{align}
for any $\varepsilon > 0$ and any $( b_0, \nu_0 ) \in \Theta$ where $B_0( 1 )$ is the $L^2$-unit ball, then weak posterior consistency holds for $\Pi$ on $\Theta$.
\end{thm}
\vskip 11pt
\begin{rmk}
In fact, consistency can be proven using the same argument we present below under the slightly weaker requirement that \eqref{kl_cond} holds for at least one representative $( b_0, \nu_0 )$ in each equivalence class of topologically equivalent points in $\Theta$.
We omit this relaxation in our proof, as it only serves to complicate terminology.
\end{rmk}
We will prove Theorem 1 by generalising the proof of Theorem 3.5 of \citep{vanderMeulen13}.
For $( b, \nu ) \in \Theta$ let $\KL( b_0, \nu_0; b, \nu )$ denote the Kullback-Leibler divergence between $p_{ \delta }^{ b_0, \nu_0 }$ and $p_{ \delta }^{ b, \nu }$:
\begin{equation*}
\KL( b_0, \nu_0; b, \nu ) := \int_{ \Omega } \int_{ \Omega } \log\left( \frac{ p_{ \delta }^{ b_0, \nu_0 }( \x, \y ) }{ p_{ \delta }^{ b, \nu }( \x, \y ) } \right) p_{ \delta }^{ b_0, \nu_0 }( \x, \y ) \pi^{ b_0, \nu_0 }( \x ) d\y d\x,
\end{equation*}
and for two probability measures $P, P'$ on the same $\sigma$-field let $K( P, P' ) := \E_P\left[ \log \left( \frac{ dP }{ dP' } \right) \right]$.
The law of a random object $Z$ under a probability measure $P$ is denoted by $\mathcal{ L }( Z | P )$.

We require the following two properties:
\begin{enumerate}
\item $\Pi( ( b, \nu ) \in \Theta : \KL( b_0, \nu_0; b, \nu ) < \varepsilon ) > 0$ for any $\varepsilon > 0$.
\item Uniform equicontinuity of the functions $\{ P_{ \delta }^{ b, \nu } f : ( b, \nu ) \in \Theta \}$ for $f \in C_b( \Omega )$, the set of bounded, continuous functions on $\Omega$.
\end{enumerate}
These two properties will be established in Lemmas \ref{kl_property} and \ref{unif_equicontinuity} below, which are the necessary generalisations of Lemmas 5.1 and A.1 of \citep{vanderMeulen13}, respectively.
\vskip 11pt
\begin{lem}\label{kl_property}
Condition \eqref{kl_cond} implies that $\Pi( ( b, \nu ) \in \Theta : \KL( b_0, \nu_0; b, \nu ) < \varepsilon ) > 0$ for any $\varepsilon > 0$.
\end{lem}
\begin{proof}
As in Lemma 5.1 of \citep{vanderMeulen13} it will be sufficient to bound $\KL( b_0, \nu_0; b, \nu )$ from above by a constant multiple of 
\begin{align*}
&\frac{ 1 }{ 2 } \E^{ b_0, \nu_0 }\left[ \| b_0( \X_0 ) - b( \X_0 ) \|_2^2 \right] + \frac{ 1 }{ 2 } \E^{ b_0, \nu_0 }\left[ \Big\| \int_{ B_0( 1 ) \setminus \{ 0 \} } \left( \frac{ d\nu_0 }{ d\nu }( \X_{ 0- }, \mathbf{ z } ) - 1 \right) \z \nu( \X_{ 0- }, d\mathbf{ z } ) \Big\|_2^2 \right]  \\
&+ \E^{ b_0, \nu_0 }\left[ \int_{ \R_0^d } \left| \log\left( \frac{ d\nu_0 }{ d\nu }( \X_{ 0- }, \mathbf{ z } ) \right) - \frac{ d \nu_0 }{ d \nu }( \X_{ 0- }, \mathbf{ z } ) + 1 \right| \nu_0( \X_{ 0- }, d\mathbf{ z } ) \right].
\end{align*}
A formal calculation yields
\begin{align}
&\int_{ \Omega } \int_{ \Omega } \log\left( \frac{ \pi^{ b_0, \nu_0 }( \mathbf{ x } ) p_{ \delta }^{ b_0, \nu_0 }( \mathbf{ x }, \mathbf{ y } ) }{ \pi^{ b, \nu }( \mathbf{ x } ) p_{ \delta }^{ b, \nu }( \mathbf{ x }, \mathbf{ y } ) } \right) p_{ \delta }^{ b_0, \nu_0 }( \mathbf{ x }, \mathbf{ y } ) \pi^{ b_0, \nu_0 }( \mathbf{ x } ) d\mathbf{ y } d\mathbf{ x } \nonumber \\
&= K( \pi^{ b_0, \nu_0 }, \pi^{ b, \nu } ) + \KL( b_0, \nu_0; b, \nu ) = K( \mathcal{ L }( \X_0, \X_{ \delta } | \mathbb{ P }^{ b_0, \nu_0 } ), \mathcal{ L }( \X_0, \X_{ \delta } | \mathbb{ P }^{ b, \nu } ) )  \nonumber \\
&\leq K( \mathcal{ L }( ( \X_t )_{ t \in [0, \delta ] } | \mathbb{ P }^{ b_0, \nu_0 } ), \mathcal{ L }( ( \X_t )_{ t \in [0, \delta ] } | \mathbb{ P }^{ b, \nu } ) ) \nonumber \\
&= K( \pi^{ b_0, \nu_0 }, \pi^{ b, \nu } ) + \E^{ b_0, \nu_0 }\left[ \log\left( \frac{ d\mathbb{ P }_{ \X_0 }^{ b_0, \nu_0 } }{ d\mathbb{ P }_{ \X_0 }^{ b, \nu } }( ( \X_t )_{ t \in [ 0, \delta ] } ) \right) \right] \label{jensen}
\end{align}
by the conditional version of Jensen's inequality.

The aim is to identify the Radon-Nikodym derivative using \citep[Theorem 2.4]{Cheridito05}, the hypotheses of which will now be verified.
The local boundedness of \citep[condition (2.1)]{Cheridito05} follow from the Lipschitz continuity assumed in \eqref{lip}.
We also require \citep[conditions (2.9) and (2.11)]{Cheridito05} --- the accompanying \citep[condition (2.10)]{Cheridito05} concerns killed processes and does not apply to our setting.
To that end, let $\{ \Omega_n \}_{ n = 1 }^{ \infty }$ denote a sequence of bounded, open subsets of $\Omega$ such that $\Omega_1 \subset \Omega_2 \subset \ldots$ and $\cup_{ n \geq 1 } \Omega_n = \Omega$.
Then
\begin{align}
&\sup_{ \mathbf{ x } \in \Omega_n } \left\{ \Big\| b_0( \mathbf{ x } ) - b( \mathbf{ x } ) - \int_{ \R_0^d } \left[ \frac{ d\nu_0 }{ d \nu }( \mathbf{ x }, \mathbf{ z } ) - 1 \right] \mathds{ 1 }_{ ( 0, 1 ] }( \| \mathbf{ z } \|_2 ) \mathbf{ z } \nu( \mathbf{ x }, d\mathbf{ z } ) \Big\|_2 \right\} \nonumber \\
&\leq \sup_{ \mathbf{ x } \in \Omega_n }\{ \| b_0( \x ) - b( \x ) \|_2 \} + \sup_{ \mathbf{ x } \in \Omega_n }\left\{ \Big\| \int_{ \R_0^d } \left[ \frac{ d\nu_0 }{ d \nu }( \mathbf{ x }, \mathbf{ z } ) - 1 \right] \mathds{ 1 }_{ ( 0, 1 ] }( \| \mathbf{ z } \|_2 ) \mathbf{ z } \nu( \mathbf{ x }, d\mathbf{ z } ) \Big\|_2 \right\} \label{cfy1}.
\end{align}
Note that the first term on the RHS is finite by Lipschitz continuity of $b_0$ and $b$, and boundedness of $\Omega_n$.
For the second, we have two cases: either $\nu( \x, \cdot )$ is a finite measure for each $\x \in \Omega$, or there exists an open set $A$ containing the origin such that $\frac{ d\nu_0 }{ d\nu }( \x, \cdot )|_A \equiv 1$ for every $\x \in \Omega$.
In the former case, Jensen's inequality gives the bound
\begin{equation*}
\sup_{ \mathbf{ x } \in \Omega_n }\left\{ \Big\| \int_{ \R_0^d } \left[ \frac{ d\nu_0 }{ d \nu }( \mathbf{ x }, \mathbf{ z } ) - 1 \right] \mathds{ 1 }_{ ( 0, 1 ] }( \| \mathbf{ z } \|_2 ) \mathbf{ z } \nu( \mathbf{ x }, d\mathbf{ z } ) \Big\|_2 \right\} \leq \sup_{ \substack{ \x \in \Omega_n, \\ \z \in B_0( 1 ) } }\left\{ \left| \frac{ d\nu_0 }{ d\nu }( \x, \z ) - 1 \right| \nu( \x, B_0( 1 ) ) \right\},
\end{equation*}
where the RHS is finite by \eqref{rn_bounds}.
If instead $\nu( \x, \cdot )$ is infinite and $\frac{ d\nu_0 }{ d\nu }( \x, \cdot )|_A \equiv 1$, then the domain of integration can be changed to $\R_0^d \setminus A$ because the integrand equals 0 on $A$.
Now $\nu( \x, \R_0^d \setminus A )$ is finite, and hence a slight modification of the previous argument shows that \eqref{cfy1} is finite in this case as well:
\begin{equation*}
\sup_{ \mathbf{ x } \in \Omega_n } \left\{ \Big\| b_0( \mathbf{ x } ) - b( \mathbf{ x } ) - \int_{ \R_0^d } \left[ \frac{ d\nu_0 }{ d \nu }( \mathbf{ x }, \mathbf{ z } ) - 1 \right] \mathds{ 1 }_{ ( 0, 1 ] }( \| \mathbf{ z } \|_2 ) \mathbf{ z } \nu( \mathbf{ x }, d\mathbf{ z } ) \Big\|_2 \right\} < K_n
\end{equation*}
for some $K_n < \infty$.
A similar argument involving the two separate cases outlined in Definition \ref{prior_supp} shows that
\begin{equation*}
\sup_{ \mathbf{ x } \in \Omega_n } \left\{ \int_{ \R_0^d } \left[ \frac{ d\nu_0 }{ d\nu }( \mathbf{ x }, \mathbf{ z } ) \log\left( \frac{ d\nu_0 }{ d\nu }( \mathbf{ x }, \mathbf{ z } ) \right) - \frac{ d\nu_0 }{ d\nu }( \mathbf{ x }, \mathbf{ z } ) + 1 \right] \nu( \mathbf{ x }, d\mathbf{ z } ) \right\} < K_n'
\end{equation*}
for some $K_n' < \infty$.
Thus, conditions (2.9) and (2.11) in \citep[Remark 2.5]{Cheridito05} hold for each $( b, \nu ) \in \Theta$, and hence \citep[Theorem 2.4]{Cheridito05} applies.

By \citep[Theorem 2.4]{Cheridito05}, the Radon-Nikodym derivative on the RHS of \eqref{jensen} can be expressed as $\E^{ b_0, \nu_0 }[ \log( \mathcal{ E }( L_{ \delta } ) ) ]$, where $\mathcal{E}$ is the Dol\'eans-Dade stochastic exponential and the process $L := ( L_t )_{ t \in [ 0, \delta ] }$ is given as
\begin{align*}
L_t =& \int_0^t \int_{ \R_0^d }\left[ \frac{ d\nu_0 }{ d\nu }( \X_{ s- }, \mathbf{ z } ) - 1 \right] ( \Z^{ \nu }( \X_{ s- }, d\mathbf{ z }, ds ) - \nu( \X_{ s- }, d\mathbf{ z } ) ds ) \\
&+ \int_0^t \left( b_0( \X_s ) - b( \X_s ) - \int_{ \R_0^d } \left( \frac{ d\nu_0 }{ d\nu }( \X_{ s- }, \mathbf{ z } ) - 1 \right) \mathds{ 1 }_{ ( 0, 1 ] }( \| \mathbf{ z } \|_2 ) \mathbf{ z } \nu( \X_{ s- }, d\mathbf{ z } ) \right) \cdot d\mathbf{ X }_s^c,
\end{align*}
where $( \X^c_s )_{ s \geq 0 }$ is the continuous martingale part of $\X$, i.e.~a $d$-dimensional Brownian motion in this setting, and $\Z^{ \nu }( \mathbf{ x }, \cdot, \cdot )$ is a Poisson random measure with intensity $\nu( \mathbf{ x }, d\mathbf{ z } ) \otimes ds$.
Note that under $\Prb^{ b_0, \nu_0 }$ the process $L$ is a local martingale, $L^c$ is a continuous local martingale with quadratic variation 
\begin{equation*}
\left< L^c \right>_t = \int_0^t \Big\| b_0( \X_s ) - b( \X_s ) - \int_{ \R_0^d } \left( \frac{ d\nu_0 }{ d\nu }( \X_{ s- }, \mathbf{ z } ) - 1 \right) \mathds{ 1 }_{ ( 0, 1 ] }( \| \mathbf{ z } \|_2 ) \mathbf{ z } \nu( \X_{ s- }, d\mathbf{ z } ) \Big\|_2^2 ds \\
\end{equation*}
and jump discontinuities of $L$ can be written as 
\begin{equation*}
\Delta L_t = \left[ \frac{ d \nu_0 }{ d \nu }( \X_{ t- }, \Delta \X_t ) - 1 \right] \mathds{ 1 }_{ ( 0, \infty ) }( \| \Delta \X_t \|_2 ),
\end{equation*}
where $\Delta \X_t$ denotes a jump discontinuity of $\X$ at time $t$.
Now, the expected quadratic variation of $L_t^c$ can be bounded by
\begin{equation*}
\E^{ b_0, \nu_0 }[ \left< L^c \right>_t ] \leq \int_0^t \E^{ b_0, \nu_0 }[ \| b_0( 0 ) + b( 0 ) + 2 C_1 \mathbf{ X }_s + K \|_2^2 ] ds
\end{equation*}
for some constant $K > 0$, using \eqref{lip}, the uniform upper and lower bounds on $\frac{ d\nu_0 }{ d\nu }$, and the fact that either $\nu$ and $\nu_0$ are equivalent and either both finite measures, or $\frac{ d \nu_0 }{ d \nu } \equiv 1$ on a neighbourhood of 0 and $\nu$ is a finite measure on any open set not containing the origin.
The stationary density has a finite second moment by Proposition \ref{existence}, so that $\E^{ b_0, \nu_0 }[ \left< L^c \right>_t ] \leq K' t$ for some other constant $K' > 0$.
Likewise, 
\begin{equation*}
\E^{ b_0, \nu_0 }\left[ \sum_{ t : \| \Delta \mathbf{ X }_t \|_2 \neq 0 } \Delta L_t^2 \right] = \int_0^t \E^{ b_0, \nu_0 }\left[ \int_{ \R_0^d } \left( \frac{ d \nu_0 }{ d \nu }( \X_{ s- }, \mathbf{ z } ) - 1 \right)^2 \nu( \X_{ s- }, d\mathbf{ z } ) \right] ds
\end{equation*}
is finite due to the aforementioned conditions on $\nu_0$ and $\nu$.
Thus $L$ has expected quadratic variation
\begin{equation*}
\E^{ b_0, \nu_0 }[ \left< L \right>_t ] = \E^{ b_0, \nu_0 }\left[ \sum_{ t : \| \Delta \mathbf{ X }_t \|_2 \neq 0 } \Delta L_t^2 + \left< L^c \right>_t \right] < \infty
\end{equation*}
for each $t > 0$, and is a true $\Prb^{ b_0, \nu_0 }$-martingale by Corollary 3 on page 73 of \citep{Protter05}.
Thus the Radon-Nikodym term in \eqref{jensen} can be written as
\begin{align}
&\E^{ b_0, \nu_0 }\Bigg[ \log\left( \frac{ d\mathbb{ P }_{ \mathbf{ x } }^{ b_0, \nu_0 } }{ d\mathbb{ P }_{ \mathbf{ x } }^{ b, \nu } }( ( \X_t )_{ t \in [ 0, \delta ] } ) \right) \Bigg] = \E^{ b_0, \nu_0 }[ \log( \mathcal{ E }( L_t ) ) ] \nonumber\\
&= \E^{ b_0, \nu_0 }\Bigg[ L_{ \delta} - L_0 - \frac{ 1 }{ 2 } \left< L^c \right>_{ \delta } + \sum_{ t : \Delta \X_t \neq 0 } \{ \log( 1 + \Delta L_{ t } ) - \Delta L_t \} \Bigg] \nonumber \\
&= \E^{ b_0, \nu_0 }\Bigg[ \frac{ -1 }{ 2 } \int_0^{ \delta } \Big\| b_0( \X_t ) - b( \X_t ) - \int_{ \R_0^d } \left( \frac{ d\nu_0 }{ d\nu }( \X_{ t- }, \mathbf{ z } ) - 1 \right) \mathds{ 1 }_{ ( 0, 1 ] }( \| \mathbf{ z } \|_2 ) \mathbf{ z } \nu( \X_{ t- }, d\mathbf{ z } ) \Big\|_2^2 dt \nonumber \\
&\;\;\;\;\;\;\;\;\;\;\;\;\;\;+ \sum_{ 0 \leq t \leq \delta : \Delta \X_t \neq 0 } \left\{ \log\left( \frac{ d\nu_0 }{ d\nu }( \X_{ t- }, \Delta \X_t ) \right) - \left( \frac{ d \nu_0 }{ d \nu }( \X_{ t- }, \Delta \X_t ) - 1 \right) \right\} \Bigg] \nonumber \\
&\leq \frac{ \delta }{ 2 } \E^{ b_0, \nu_0 }\left[ \| b_0( \X_0 ) - b( \X_0 ) \|_2^2 \right] + \frac{ \delta }{ 2 } \E^{ b_0, \nu_0 }\left[ \Big\| \int_{ B_0( 1 ) \setminus \{ 0 \} } \left( \frac{ d\nu_0 }{ d\nu }( \X_{ 0- }, \mathbf{ z } ) - 1 \right) \z \nu( \X_{ 0- }, d\mathbf{ z } ) \Big\|_2^2 \right]  \nonumber \\
&\;\;\;\;\;\;\;+ \delta \E^{ b_0, \nu_0 }\left[ \int_{ \R_0^d } \left| \log\left( \frac{ d\nu_0 }{ d\nu }( \X_{ 0- }, \mathbf{ z } ) \right) - \frac{ d \nu_0 }{ d \nu }( \X_{ 0- }, \mathbf{ z } ) + 1 \right| \nu_0( \X_{ 0- }, d\mathbf{ z } ) \right], \label{kl_bound}
\end{align}
where the first equality follows from Theorem 2.4 of \citep{Cheridito05}, the second by definition of $\mathcal{ E }$ for jump diffusion processes, and the remainder of the calculation by stationarity and because $\nu_0$ is the compensator of the Poisson random measure driving the jumps of $\X$ under $\mathbb{ P }^{ b_0, \nu_0 }$.
The result now follows from \eqref{jensen} and \eqref{kl_bound}.
\end{proof}
\begin{lem}\label{unif_equicontinuity}
For each $\delta > 0$ and $f \in C_b( \Omega )$, the collection $\{ P_{ \delta }^{ b, \nu } f : ( b, \nu ) \in \Theta \}$ is locally uniformly equicontinuous: for any compact $K \in \Omega$ and $\varepsilon > 0$ there exists $\gamma := \gamma( \varepsilon, K, f, \delta ) > 0$ such that 
\begin{equation*}
\sup_{ ( b, \nu ) \in \Theta } \sup_{ \substack{\x, \y \in K : \\ \| \x - \y \|_2 < \gamma } } | P_{ \delta }^{ b, \nu } f( \x ) - P_{ \delta }^{ b, \nu } f( \y ) | < \varepsilon.
\end{equation*}
\end{lem}
\begin{proof}
Theorem 2.3 of \citep{Wang10} establishes Lipschitz continuity for jump diffusions satisfying \eqref{lip} using a coupling argument for $f \in \mathcal{ B }_b( \Omega )$, the set of bounded, measurable functions.
We begin by showing that the conditions of \citet{Wang10} are satisfied.

In our notation and setting, the condition of Theorem 2.3 of \citep{Wang10} is that for some constant $\beta \in ( 0, 1 )$ there exists a constant $C_{ \beta } > 0$ such that 
\begin{align*}
&\frac{ ( 1 + \| \x - \y \|_2 ) \left< b( \x ) - b( \y ), \x - \y \right> }{ \|\x - \y \|_2 } + \frac{ ( 1 + \| \x - \y \|_2 ) \int_{ \| \z \|_2 \leq 1 } \| c( \x, \z ) - c( \y, \z ) \|_2^2 M( d\z ) }{ 2 \| \x - \y \|_2 } \\
&+ ( 1 + \| \x - \y \|_2 ) \int_{ \| \z \|_2 > 1 } \| c( \x, \z ) - c( \y, \z ) \|_2 M( d\z ) + ( 1 + \| \x - \y \|_2 ) C_{ \beta } \leq 2
\end{align*}
whenever $\| \mathbf{ x } - \mathbf{ y } \|_2 < \beta$ and where $\left< \cdot, \cdot \right>$ denotes the usual Euclidean inner product. 
By \eqref{lip}, the first two summands on the LHS can be bounded by $\beta ( 1 + \beta ) \sqrt{ C_1 }$ and $\beta ( 1 + \beta ) C_1 / 2$, respectively.
The fourth is trivially bounded by $( 1 + \beta ) C_{ \beta }$.
By Jensen's inequality, \eqref{lip} and \eqref{m_mass}, the third term can be bounded by $\beta ( 1 + \beta ) \sqrt{ C_1 C_5 }$.
Hence, the whole LHS can be bounded by
\begin{equation*}
( 1 + \beta ) \left[ \beta \sqrt{ C_1 } \left( 1 + \sqrt{ C_5 } \right) + \frac{ \beta C_1 }{ 2 } + C_{ \beta } \right]
\end{equation*}
which can clearly be made arbitrarily small by choosing both $\beta$ and $C_{ \beta }$ to be sufficiently small.
This choice can be made uniformly due to the uniform bounds on the Lipschitz constant $C_1$ and the total mass constraint $C_5$.

Now, the Lipschitz constant in \citep[Theorem 2.3]{Wang10} is of the form
\begin{equation*}
2 ( 1 + \beta ) \left( \frac{ 1 }{ C_{ \beta } \Delta } + \frac{ 1 }{ \beta } \right) \| f \|.
\end{equation*}
Since $\beta$ and $C_{ \beta }$ can be chosen uniformly in $\Theta$, and $\Delta$ and $\| f \|$ are constants, uniform equicontinuity holds.
\end{proof}
\begin{proof}[Proof of Theorem \ref{ts_consistency}]
Given Lemmas \ref{kl_property} and \ref{unif_equicontinuity} above, the remainder of the proof follows as in \citep{vanderMeulen13}.
It suffices to show that for $f \in C_b( \Omega )$ and $B := \{ ( b, \nu ) \in \Theta : \| P_{ \delta }^{ b, \nu } f - P_{ \delta }^{ b_0, \nu_0 } f \|_{ 1, \rho } > \varepsilon \}$ we have $\Pi( B | \x_{ 0 : n } ) \rightarrow 0$ with $\Prb^{ b_0, \nu_0 }$-probability 1.
To that end we fix $f \in \operatorname{Lip}( \Omega )$ and $\varepsilon > 0$ and thus the set $B$.
Lemma \ref{kl_property} implies that Lemma 5.2 of \citep{vanderMeulen13} holds, so that if, for measurable subsets $C_n \subset \Theta$, there exists $c > 0$ such that
\begin{equation*}
e^{ n c } \int_{ C_n } \pi^{ b, \nu }( \x_0 ) \prod_{ i = 1 }^n p_{ \delta }^{ b, \nu }( \x_{ i - 1 }, \x_i ) \Pi( db, d\nu ) \rightarrow 0
\end{equation*}
$\Prb^{ b_0, \nu_0 }$-a.s.~then $\Pi( C_n | \x_{ 0 : n } ) \rightarrow 0$ $\Prb^{ b_0, \nu_0 }$-a.s.~as well.
Likewise, Lemma \ref{unif_equicontinuity} implies Lemma 5.3 of \citep{vanderMeulen13}: there exists a compact subset $K \subset \Omega$, $N \in \N$ and compact, connected sets $I_1, \ldots, I_N$ that cover $K$ such that 
\begin{equation*}
B \subset \bigcup_{ j = 1 }^N B_j^+ \cup \bigcup_{ j = 1 }^N B_j^-,
\end{equation*}
where
\begin{align*}
B_j^+ &:= \left\{ ( b, \nu ) \in \Theta : P_{ \delta }^{ b, \nu } f( \x ) - P_{ \delta }^{ b_0, \nu_0 } f( \x ) > \frac{ \varepsilon }{ 4 \rho( K ) } \text{ for every } \x \in I_j \right\}, \\
B_j^- &:= \left\{ ( b, \nu ) \in \Theta : P_{ \delta }^{ b, \nu } f( \x ) - P_{ \delta }^{ b_0, \nu_0 } f( \x ) < \frac{ -\varepsilon }{ 4 \rho( K ) } \text{ for every } \x \in I_j \right\}.
\end{align*}
Thus it is only necessary to show $\Pi( B_j^{ \pm } | \x_{ 0 : n } ) \rightarrow 0$ $\Prb^{ b_0, \nu_0 }$-almost surely.
Define the stochastic process
\begin{equation*}
D_n := \left( \int_{ B_j^+ } \pi^{ b, \nu }( \x_0 ) \prod_{ i = 1 }^n p_{ \delta }^{ b, \nu }( \x_{ i - 1 }, \x_i ) \Pi( db, d\nu ) \right)^{ 1 / 2 }.
\end{equation*}
Now $D_n \rightarrow 0$ exponentially fast as $n \rightarrow \infty$ by an argument identical to that used to prove Theorem 3.5 of \citep{vanderMeulen13}.
The same is also true of the analogous stochastic process defined by integrating over $B_j^-$, which completes the proof.
\end{proof}

We conclude this section with a corollary giving explicit sufficient conditions for the hypotheses of Theorem \ref{ts_consistency}.
\vskip 11pt
\begin{cor}\label{explicit_conds}
Suppose the prior is of the product form $\Pi = \Pi_b \otimes \Pi_c \otimes \Pi_M$, with the three factors corresponding to the drift $b$, the jump function $c$ and the jump measure $M$ in the obvious way, and let their respective supports be $\Theta_b$, $\Theta_c$, and $\Theta_M$.
Then the following numerous but simple conditions are sufficient for posterior consistency:
\begin{enumerate}
\item There exists a $K_1 > 0$ such that $\| b( \x ) - b( \y ) \|_2^2 \leq K_1 \| \x - \y \|_2^2$ uniformly in $b \in \Theta_b$. \label{cor_1}
\item For each $b \in \Theta_b$ there exists a $K_b > 0$ and an $R_b > 0$ such that $\| \x \|_2 > R_b \Rightarrow \x \cdot b( \x ) \leq - K_b \| \x \|_2^2$. \label{cor_2}
\item For any compact set $A \subset \Omega$ and any $b_0 \in \Theta_b$, $\Pi_b( b \in \Theta_b : \sup_{ \x \in A } \left\{ \| b( \x ) - b_0( \x ) \|_2^2 \right\} < \varepsilon ) > 0$ for any $\varepsilon > 0$. \label{cor_3}
\item There exists a $K_2 > 0$ such that $\| c( \x, \z ) - c( \y, \z ) \|_2^2 \leq K_2 \| \x - \y \|_2^2$ for every $\z \in \R_0^d$, uniformly in $c \in \Theta_c$. \label{cor_5}
\item For every $c \in \Theta_c$ there exists a $K_c > 0$ such that $\| c( \x, \z ) - c( \x, \bxi ) \|_2^2 \leq K_c \| \z - \bxi \|_2^2$ for every $\x \in \Omega$. \label{cor_6}
\item For every $c \in \Theta_c$ we have $c( \cdot, 0 ) \equiv 0$. \label{cor_5.5}
\item For every $\x \in \Omega$ there exists a compact set $J_{ \x } \subset \R_0^d$ which is the image of $c( \x, \cdot )$ simultaneously for every $c \in \Theta_c$. \label{cor_7}
\item For each $c \in \Theta_c$ there exists an $R_c > 0$ such that 
\begin{equation*}
\sup_{ \x \in \Omega, \z \in \R_0^d }\{ c( \x, \z ) \} = \sup_{ \x \in B_0( R_c ), \z \in \R_0^d }\{ c( \x, \z ) \},
\end{equation*}
and likewise for infima. \label{cor_new}
\item $c( \mathbf{ x }, \cdot )$ is not constant on any open ball in $\R_0^d$ for any $\mathbf{ x } \in \Omega$ or any $c \in \Theta_c$. \label{cor_8}
\item For any compact set $A \subset \Omega \times \R_0^d$ and any $c_0 \in \Theta_c$, 
\begin{equation*}
\Pi_c( c \in \Theta_c : \sup_{ ( \x, \z ) \in A }\left\{ \| c( \x, \z ) - c_0( \x, \z ) \|_2 \right\} < \varepsilon ) > 0
\end{equation*}
for any $\varepsilon > 0$. \label{cor_9}
\item There exists a compact set $J \subset \R_0^d$ such that every $M \in \Theta_M$ has $M( J^c ) = 0$, and there exists $K_3 < \infty$ such that $M( J ) \leq K_3$ uniformly in $M \in \Theta_M$. \label{cor_11}
\item Each $M \in \Theta_M$ has a density $M( \z )$ on $J$ that is bounded above and away from 0. \label{cor_13}
\item For any $M_0 \in \Theta_M$, $\Pi_M( M : \| M - M_0 \|_{ \infty } < \varepsilon ) > 0$ for any $\varepsilon > 0$. \label{cor_14}
\end{enumerate}
\end{cor}
\vskip 11pt
\begin{rmk}
Conditions \ref{cor_5} -- \ref{cor_8} characterising functions $c \in \Theta_c$ are perhaps the least transparent of the three sets in the above corollary.
A straightforward, one-dimensional example of a function which satisfies them can be built from
\begin{equation*}
c( x, z ) = 2 + \sin( z ),
\end{equation*}
by using any Lipschitz modification to ensure that $c(x, 0) = 0$.
More generally, families of functions with no dependence on the first argument, with bounded gradient in the second argument, with no intervals of constancy, and which map into a compact set will satisfy these conditions provided that they are thus modified at 0.
\end{rmk}
\begin{proof}[Proof of Corollary \ref{explicit_conds}]
We begin by verifying the hypotheses of Proposition \ref{existence}.
By conditions \ref{cor_1}, \ref{cor_5}, and \ref{cor_11}, the bound
\begin{align*}
\| b( \mathbf{ x } ) - b( \mathbf{ y } ) \|_2^2 + \int_J \| c( \mathbf{ x }, \z ) - c( \mathbf{ y }, \z ) \|_2^2 M( d\z ) &\leq K_1 \| \x - \y \|_2^2 + K_2 \int_J \| \x - \y \|_2^2 M( d \z ) \\
&\leq ( K_1 + K_2 K_3 ) \| \x - \y \|_2^2
\end{align*}
holds uniformly in $b \in \Theta_b$, $c \in \Theta_c$ and $M \in \Theta_M$.
Thus, \eqref{lip} holds.
Furthermore, \eqref{lip2} and \eqref{drift}, coincide with conditions \ref{cor_2} and \ref{cor_6}, and \eqref{m_moment} is implied by \ref{cor_11}.
The requirement that $c( \cdot, 0 ) \equiv 0$ holds by condition \ref{cor_5.5}.

Next, we will check that the conditions of Definition \ref{prior_supp} are satisfied.
Condition \ref{cor_11} implies \eqref{m_mass}.
To see that any two measures $\nu( \x, \cdot )$ and $\nu_0( \x, \cdot )$ are equivalent, note that neither measure can have atoms by conditions \ref{cor_8} and \ref{cor_13}, and that the support of both measures coincides by conditions \ref{cor_7} and \ref{cor_13}.
Every $\nu( \x, \cdot )$ is also finite by conditions \ref{cor_7} and \ref{cor_11}.
To see \eqref{rn_bounds}, note that by condition \ref{cor_new} we have
\begin{equation*}
\inf_{ \substack{ \x \in \Omega \\ \z \in J } } \left\{ \frac{ d\nu_0 }{ d\nu }( \x, \z ) \right\} = \inf_{ \substack{ \x \in B_0( R_c \vee R_{ c_0 } ) \\ \z \in J } } \left\{ \frac{ M_0( \{ \y \in J : c_0( \x, \y ) = \z \} ) }{ M( \{ \y \in J : c( \x, \y ) = \z \} ) } \right\}.
\end{equation*}
Both the numerator and denominator are strictly positive for every $( \x, \z ) \in B_0( R_c \vee R_{ c_0 } ) \times J$ by conditions \ref{cor_7} and \ref{cor_11}, and hence so is the infimum by compactness of $B_0( R_c \vee R_{ c_0 } ) \times J$.
The supremum is also finite by an analogous argument, verifying \eqref{rn_bounds}.

It only remains to verify \eqref{kl_cond}, which we will do by considering the three expectations separately.
Firstly,
\begin{align*}
\E^{ b_0, \nu_0 }&\left[ \| b_0( \X_0 ) - b( \X_0 ) \|_2^2 \right] = \int_{ \Omega } \| b_0( \x ) - b( \x ) \|_2^2 \pi^{ b_0, \nu_0 }( \x ) d\x \\
&\leq \int_{ B_0( n ) } \| b_0( \x ) - b( \x ) \|_2^2 \pi^{ b_0, \nu_0 }( \x ) d\x + \int_{ B_0( n )^c } \| b_0( \x ) - b( \x ) \|_2^2 \pi^{ b_0, \nu_0 }( \x ) d\x \\
&\leq \sup_{ \x \in B_0( n ) }\left\{ \| b( \x ) - b_0( \x ) \|_2^2 \right\} + \int_{ B_0( n )^c } \left\{ \| b_0( 0 ) \|_2^2 + \| b( 0 ) \|_2^2 + 2 C_1 \| \x \|_2^2 \right\} \pi^{ b_0, \nu_0 }( \x ) d\x.
\end{align*}
The second term on the RHS vanishes as $n \rightarrow 0$ because $\pi^{ b_0, \nu_0 }$ has a second moment by Proposition \ref{existence}.
Thus, a positive probability of arbitrarily small values of $\sup_{ \x \in B_0( n ) }\{ \| b( \x ) - b_0( \x ) \|_2^2 \}$ for any finite $n$, assumed in condition \ref{cor_3}, implies a positive probability of arbitrarily small values of $\E^{ b_0, \nu_0 }[ \| b_0( \X_0 ) - b( \X_0 ) \|_2^2 ]$.

Next, Jensen's inequality and condition \ref{cor_11} yields
\begin{align*}
&\E^{ b_0, \nu_0 }\left[ \Big\| \int_{ B_0( 1 ) \cap J } \left( \frac{ d\nu_0 }{ d\nu }( \X_{ 0- }, \mathbf{ z } ) - 1 \right) \z \nu( \X_{ 0- }, d\mathbf{ z } ) \Big\|_2^2 \right] \\
&\leq \E^{ b_0, \nu_0 }\left[ \nu( \X_{ 0- }, B_0( 1 ) \cap J ) \int_{ B_0( 1 ) \cap J } \left( \frac{ d\nu_0 }{ d\nu }( \X_{ 0- }, \z ) - 1 \right)^2 \nu( \X_{ 0- }, d\z ) \right] \\
&\leq \sup_{ \substack{ \x \in B_0( R_c \vee R_{ c_0 } ) \\ \z \in B_0( 1 ) \cap J } } \left\{ \nu( \x, B_0( 1 ) \cap J )^2 \left( \frac{ d\nu_0 }{ d\nu }( \x, \z ) - 1 \right)^2 \right\} \leq K \sup_{ \substack{ \x \in B_0( R_c \vee R_{ c_0 } ) \\ \z \in B_0( 1 ) \cap J } } \left\{ \left( \frac{ d\nu_0 }{ d\nu }( \x, \z ) - 1 \right)^2 \right\}
\end{align*}
for some $K < \infty$, as $\sup_{ \x \in B_0( R_c \vee R_{ c_0 } ) }\left\{ \nu( \x, B_0( 1 ) \cap J )^2 \right\}$ is finite because $M$ is finite, $c$ is Lipschitz continuous, and $B_0( R_c \vee R_{ c_0 } )$ is compact.
Now, for a measurable set $A \in \mathcal{ B }( J )$ and $\x \in B_0( R_c \vee R_{ c_0 } )$, we have that
\begin{equation*}
\frac{ d\nu_0 }{ d\nu }( \x, A ) = \frac{ M_0( \{ \y \in J : c_0( \x, \y ) \in A \} ) }{ M( \{ \y \in J : c( \x, \y ) \in A \} ) } \leq \frac{ M( \{ \y \in J : c_0( \x, \y ) \in A \} ) + | J | \| M - M_0 \|_{ \infty } }{ M( \{ \y \in J : c( \x, \y ) \in A \} ) }.
\end{equation*}
Now, for each point $\y \in J : c_0( \x, \y ) \in A$ we also have $d( c( \x, \y ), A ) \leq \| c_0( \x, \y ) - c( \x, \y ) \|_2$, where $d( \x, A ) := \inf_{ \y \in A } \| \x - \y \|_2$.
Thus, if 
\begin{equation*}
\sup_{ \x \in B_0( R_c \vee R_{ c_0 } ), \y \in J }\{ \| c_0( \x, \y ) - c( \x, \y ) \|_2 \} < \varepsilon,
\end{equation*}
then 
\begin{equation*}
\{ \y \in J : c_0( \x, \y ) \in A \} \subset \{ \y \in J : d( c( \x, \y ), A ) \leq \varepsilon \}.
\end{equation*}
By continuity, the $M$-mass of the latter set converges to $M( \{ \y \in J : c( \x, \y ) \in A \} )$ as $\varepsilon \rightarrow 0$, so that $( \frac{ d\nu_0 }{ d\nu }( \x, \z ) - 1 )^2$ can be made arbitrarily small uniformly in $( \x, \z ) \in  B_0( R_c \vee R_{ c_0 } ) \times (J \cap B_0( 1 ) )$ by choosing $c$ such that 
\begin{equation*}
\sup_{ \x \in B_0( R_c \vee R_{ c_0 } ), \z \in J \cap B_0( 1 ) }\left\{ \| c( \x, \z ) - c_0( \x, \z ) \|_2 \right\} < \varepsilon,
\end{equation*}
and $M$ such that $ \| M - M_0 \|_{ \infty } < \varepsilon'$, for sufficiently small constants $\varepsilon > 0$, and $\varepsilon' > 0$.
By conditions \ref{cor_9} and \ref{cor_14}, the prior assigns positive probability to both of these events, and as they are independent, the joint probability is positive as well.

The proof that the third expectation in \eqref{kl_cond},
\begin{equation*}
\E^{ b_0, \nu_0 }\left[ \int_J \left| \log\left( \frac{ d\nu_0 }{ d\nu }( \X_{ 0- }, \mathbf{ z } ) \right) - \frac{ d \nu_0 }{ d \nu }( \X_{ 0- }, \mathbf{ z } ) + 1 \right| \nu_0( \X_{ 0- }, d\mathbf{ z } ) \right],
\end{equation*}
takes arbitrarily small values with strictly positive probability under $\Pi$ follows from a calculation very similar to the second one above, and is omitted.
\end{proof}

\section{Example priors}\label{prior}

The conditions of Theorem \ref{ts_consistency} and Corollary \ref{explicit_conds} are verifiable, in that they do not depend on intractable quantities, but it is not immediately clear whether there exists a prior $\Pi$ that satisfies them.
In this section we show that there are at least three distributions which can be used to construct prior of the form considered in Corollary \ref{explicit_conds}:  the wavelet prior \citep{Ruggeri05}, the discrete net prior \citep{Ghosal97}, and the Dirichlet mixture model prior \citep{Lo84}.
We will show that these three priors satisfy the respective conditions on $\Pi_b$, $\Pi_c$, and $\Pi_M$ in Corollary \ref{explicit_conds}, though it will also be clear that either of the first two priors could be used to construct any of the three components.
The Dirichlet process mixture model prior can only used as the $\Pi_M$ component without additional assumptions as it does not allow for control of the growth rate of samples.

Discrete net priors were also used by both \citet{vanderMeulen13}, and \citet{Gugushvili14} to demonstrate the existence of consistent priors for nonparametric inference of drifts for diffusions.
Wavelet priors were also used in the one-dimensional setting in \citep{vanderMeulen13}, and our calculation demonstrates that they remain tractable in any dimension.

\subsection{Wavelet priors}\label{wavelet_prior}

Let $\bbZ_{ \Omega }^d := \bbZ^d \cap \Omega$, and for each $j \geq 0$ define the sets 
\begin{equation*}
\Lambda_j := \{ 2^{ -j } \bfk + 2^{ - j - 1 } \bfzeta : \bfk \in \bbZ_{ \Omega }^d, \bfzeta \in \{ 0, 1 \}^d \setminus \{ 0, \ldots, 0 \} \}.
\end{equation*}
Let 
\begin{equation*}
( \{ \phi( \cdot  - \bfk ) \}_{ \bfk \in \bbZ_{ \Omega }^d }, \{ 2^{ j d / 2 } \psi( 2^j \cdot - \bflambda ) \}_{ \bflambda \in \Lambda_j, j \geq 0 } )
\end{equation*}
 be an orthonormal wavelet basis of $L^2( \Omega )$, in other words be such that every $f \in L^2( \Omega )$ can be written as
\begin{equation*}
f( \x ) = \sum_{ \bfk \in \bbZ_{ \Omega }^d } a_{ \bfk } \phi( \x - \bfk ) + \sum_{ j = 0 }^{ \infty } \sum_{ \bflambda \in \Lambda_j } a_{ \bflambda }'  2^{ j d / 2 }\psi( 2^j \x - \bflambda )
\end{equation*}
for some sequences of coefficients $\{ a_{ \bfk } \}$ and $\{ a_{ \bflambda }' \}$.
Suppose further that the wavelets $\phi$ and $\psi$ are continuously differentiable and compactly supported, and hence also that both the wavelets and their gradients are bounded functions.

Under regularity conditions on the wavelet basis \citep[Definition 2, page 21]{Meyer95}, it is well known \citep[Theorem 5, page 179]{Meyer95} that $f \in C^s$, the space of $s$-H\"older continuous functions, for $s > 0$ if and only if there exist constants $C_f > 0$ and $K_f > 0$ such that
\begin{align*}
| a_{ \bfk } | &\leq C_f \text{ for } \bfk \in \bbZ_{ \Omega }^d \\
| a_{ \bflambda }' | &\leq 2^{ - j ( d / 2 + s ) } K_f \text{ for } \bflambda \in \Lambda_j \text{ and } j \geq 0.
\end{align*}
For such a function, we have
\begin{align*}
| \nabla \cdot f( \x ) | &\leq \sum_{ \bfk \in \bbZ_{ \Omega }^d } | a_{ \bfk } | | \nabla \cdot \phi( \x - \bfk ) | + \sum_{ j = 0 }^{ \infty } \sum_{ \bflambda \in \Lambda_j } | a_{ \bflambda }' | 2^{ j d / 2 + 1 } | ( \nabla \cdot \psi )|_{ 2^j \x - \bflambda } | \\
&\leq C_f \sum_{ \bfk \in \bbZ_{ \Omega }^d } | \nabla \cdot \phi( \x - \bfk ) | + K_f \sum_{ j = 0 }^{ \infty } 2^{ - j ( s - 1 ) } \sum_{ \bflambda \in \Lambda_j } | ( \nabla \cdot \psi )|_{ 2^j \x - \bflambda } |.
\end{align*}
Define $R := \inf\{ r > 0 : ( \operatorname{supp}( \phi ) \cup \operatorname{supp}( \psi ) ) \subseteq B_0( r ) \}$, which is finite by the assumed compact support of the wavelets.
Then, for each fixed $\x \in \Omega$ we have the bounds
\begin{align*}
\# \{ \bfk \in \bbZ_{ \Omega }^d : \phi( \x - \bfk )  \neq 0 \} &\leq ( 2 R )^d, \\
\# \{ \bflambda \in \Lambda_j : \psi( 2^j \x - \bflambda )  \neq 0 \} &\leq ( 2^d - 1 ) ( 2^{ j + 1 } R )^d.
\end{align*}
These bounds, along with compact support and continuous differentiability of wavelets also give that
\begin{equation*}
| \nabla \cdot f( \x ) | \leq C_f ( 2 R )^d \| \nabla \cdot \phi \|_{ \infty } + K_f  ( 2^d - 1 ) ( 2 R )^d \| \nabla \cdot \psi \|_{ \infty } \sum_{ j = 0 }^{ \infty } 2^{ j ( d + 1 - s ) }.
\end{equation*}
which is finite provided $s > d + 1$.
Thus, for any $\gamma > s - d - 1 > 0$ and $L > 0$, the set
\begin{align*}
\mathcal{ F }_{ \gamma, L } := \Big\{ f \in L^2( \Omega ) : f( \x ) &= \sum_{ \bfk \in \bbZ_{ \Omega }^d } a_{ \bfk } \phi( \x - \bfk ) + \sum_{ j = 0 }^{ \infty } \sum_{ \bflambda \in \Lambda_j } a_{ \bflambda }'  2^{ j d / 2 }\psi( 2^j \x - \bflambda ), \\
&\sup_{ \bfk \in \bbZ_{ \Omega }^d } | a_{ \bfk } | + \sup_{ j \geq 0 } \sup_{ \bflambda \in \Lambda_j } 2^{ j ( 3 d / 2 + 1 + \gamma ) }| a_{ \bflambda }' | < L \Big\}
\end{align*}
consists of functions $f \in L^2( \Omega )$ with uniformly bounded gradient.

We specify a prior $\Pi'$ on $\mathcal{ F }_{ \gamma, L }$ as the law of the random function
\begin{equation*}
\x \mapsto \sum_{ \bfk \in \bbZ_{ \Omega }^d } U_{ \bfk } \phi( \x - \bfk ) + \sum_{ j = 0 }^{ \infty } \sum_{ \bflambda \in \Lambda_j } U_{ \bflambda }' 2^{ - j ( d + 1 + \gamma ) } \psi( 2^j \x - \bflambda ),
\end{equation*}
where $\{ U_{ \bfk } \}_{ \bfk \in \bbZ_{ \Omega }^d }$ and $\{ U_{ \bflambda }' \}_{ \bflambda \in \Lambda_j, j \geq 0 }$ are all independent random variables with uniform distributions on $( - L, L )$.
A prior on functions taking values in $\R^d$ is specified by drawing each of the $d$ coordinate maps independently from $\Pi'$.
We denote the resulting law by $( \Pi' )^{ \otimes d }$.

Now let $\{ p_m \}_{ m \in \N }$ be a strictly positive probability mass function independent of $( \Pi' )^{ \otimes d }$.
A function $b \sim \Pi_b$ is drawn by first sampling $b' \sim (\Pi')^{ \otimes d }$ and $m \sim p_m$, and setting
\begin{equation*}
b( \x ) | b', m = \begin{cases}
b'( \x ) &\text{ if } \| \x \|_2 \leq m, \\
b'( P_m \x ) - C \x  &\text{ if } \| \x \|_2 > m, \\
\end{cases}
\end{equation*}
where $P_m$ denotes the orthogonal projection to the closed ball $B_0( m )$, and $C > 0$ is a fixed constant.
Samples from $\Pi_b$ satisfy conditions \ref{cor_1} and \ref{cor_2} of Corollary \ref{explicit_conds} with probability 1 by construction.
Thus, it only remains to verify condition \ref{cor_3}.

Fix a compact $A \subset \Omega$ and a drift function $b_0 \in ( \mathcal{ F }_{ \gamma, L } )^d$ satsfying conditions \ref{cor_1} and \ref{cor_2} of Corollary \ref{explicit_conds}.
Let $m \in \N$ be such that $A \subseteq B_0( m )$, and let $b \sim \Pi_b$.
Suppose $b_{ 0, i }$, the $i^{\text{th}}$ coordinate of $b_0$, has wavelet coefficients $\{ a_{ \bfk, i }^0 \}_{ \bfk \in \bbZ_{ \Omega }^d }$ and $\{ ( a_{ \bflambda, i }^0 )' \}_{ \bflambda \in \Lambda_j, j \geq 0 } \equiv \{ 2^{ - j ( d + 1 + \gamma ) } ( U_{ \bflambda, i }^0 )' \}_{ \bflambda \in \Lambda_j, j \geq 0 }$, and likewise that random wavelet coefficients of $b_i$ are denoted by $\{ a_{ \bfk, i } \}_{ \bfk \in \bbZ_{ \Omega }^d }$ and $\{ a_{ \bflambda, i }' \}_{ \bflambda \in \Lambda_j, j \geq 0 } \equiv \{ 2^{ - j ( d + 1 + \gamma ) } U_{ \bflambda, i }' \}_{ \bflambda \in \Lambda_j, j \geq 0 }$.
Then
\begin{align}
&\sup_{ \x \in A } \{ \| b( \x ) - b_0( \x ) \|_2^2 \} \leq d \sup_{ \x \in A } \{ \| b( \x ) - b_0( \x ) \|_{ \infty }^2 \} \nonumber \\
&\leq d \sup_{ \x \in A } \Bigg\{ \max_{ i \in \{ 1, \ldots, d \} } \Bigg\{ \sum_{ \bfk \in \bbZ_{ \Omega }^d } | a_{ \bfk, i } - a_{ \bfk, i }^0 | | \phi( \x - \bfk ) | + \sum_{ j = 0 }^{ \infty } \sum_{ \bflambda \in \Lambda_j } 2^{ j d / 2 }| a_{ \bflambda, i }' - ( a_{ \bflambda, i }^0 )' | | \psi( 2^j \x - \bflambda ) | \Bigg\} \Bigg\} \nonumber \\
&\leq d \sup_{ \x \in A } \Bigg\{ \max_{ i \in \{ 1, \ldots, d \} } \Bigg\{ ( 2 R )^d \| \phi \|_{ \infty } \max_{ \bfk \in B_{ \x }( R ) \cap \bbZ_{ \Omega }^d } \{ | a_{ \bfk, i } - a_{ \bfk, i }^0 | \} \nonumber \\
&\hskip 100pt + ( 2 R )^d ( 2^d - 1 ) \| \psi \|_{ \infty } \sum_{ j = 0 }^{ \infty } 2^{ - j ( 1+ \gamma ) } \max_{ \bflambda \in B_{ \x }( R ) \cap \Lambda_j }\{ | U_{ \bflambda, i }' - ( U_{ \bflambda, i }^0 )' | \Bigg\} \Bigg\} \nonumber \\
&\leq d \sup_{ \x \in A } \Bigg\{ \max_{ i \in \{ 1, \ldots, d \} } \Bigg\{ ( 2 R )^d \| \phi \|_{ \infty } \max_{ \bfk \in B_{ \x }( R ) \cap \bbZ_{ \Omega }^d } \{ | a_{ \bfk, i } - a_{ \bfk, i }^0 | \} \nonumber \\
&\hskip 60pt + 2 ( 2 R )^d ( 2^d - 1 ) \| \psi \|_{ \infty } \left( \max_{ \bflambda \in B_{ \x }( R ) \cap ( \Lambda_0 \cup \cdots \cup \Lambda_{ j^* } ) }\{ | U_{ \bflambda, i }' - ( U_{ \bflambda, i }^0 )' | \} + L 2^{ - j^* } \right) \Bigg\} \Bigg\}, \label{cor_1_bound}
 \end{align}
 for any $j^* \in \N$, where the last inequality follows by splitting the summation over $j \in \N$ into $j \leq j^*$ and $j > j^*$.
 Now, the random variables
 \begin{align*}
&\sup_{ \x \in A } \left\{ \max_{ i \in \{ 1, \ldots, d \} } \left\{ \max_{ \bfk \in B_{ \x }( R ) \cap \bbZ_{ \Omega }^d } \{ | a_{ \bfk, i } - a_{ \bfk, i }^0 | \} \right\} \right\}, \text{ and }\\
& \sup_{ \x \in A } \left\{ \max_{ i \in \{ 1, \ldots, d \} } \left\{ \max_{ \bflambda \in B_{ \x }( R ) \cap ( \Lambda_0 \cup \cdots \cup \Lambda_{ j^* } ) }\{ | U_{ \bflambda, i }' - ( U_{ \bflambda, i }^0 )' | \}  \right\} \right\}
 \end{align*}
each take arbitrarily small values with positive probability for every $j^* \in \N$, because the suprema and maxima are taken over compact sets and thus the number of involved i.i.d.~$U( - L, L )$-distributed random variables is finite.
Thus, for each $\varepsilon > 0$ we can choose a $j^* \in \N$ such that the RHS of \eqref{cor_1_bound} is smaller than $\varepsilon$ with positive probability, and hence condition \ref{cor_3} of Corollary \ref{explicit_conds} holds.

\subsection{Discrete net priors}\label{net_prior}

Let $J \subset \mathbb { R }_0^d$ be a compact set excluding the origin, and let $\Theta_c'$ be a set of functions $c : \Omega \times J \mapsto J$ satisfying conditions \ref{cor_5}, \ref{cor_6}, \ref{cor_7}, and \ref{cor_8} of Corollary \ref{explicit_conds}, and with the constant $K_c$ in condition \ref{cor_6} also holding uniformly in $\Theta_c'$.
Let $\Theta_c^{ ( m ) } := \{ c|_{ B_0( m ) } : c \in \Theta_c' \}$ be the set of restrictions of the first coordinate to the closed ball $B_0( m ) \subset \Omega$.
By uniform equicontinuity and the Arzel\`a-Ascoli theorem, $\Theta_c^{ ( m ) }$ is totally bounded in the uniform norm.
Hence, for every $n$, it is possible to construct a finite $\varepsilon_n$-net $\Theta_c^{ ( m, n ) }$ over $\Theta_c^{ ( m ) }$, where $\{ \varepsilon_n \}_{ n \in \mathbb{ N } }$ is a sequence of strictly positive numbers tending to 0.
In other words, $\Theta_c^{ ( m, n ) }$ is a finite set with the property that every element of $\Theta_c^{ ( m ) }$ is within distance $\varepsilon_n$ of some element of $\Theta_c^{ ( m, n ) }$ in the supremum norm.
Each element of $\Theta_c^{ ( m, n ) }$ can then be extended to a function on the whole $\Omega \times J$ by setting $c( \mathbf{ x }, \mathbf{ z } ) = c( P_m \mathbf{ x }, \mathbf{ z } )$ outside $B_0( m )$, where $P_m$ again denotes orthogonal projection to $B_0( m )$.
A discrete net prior is constructed by fixing two strictly positive probability mass functions, $\{ p_m \}_{ m \in \mathbb{ N } }$ and $\{ q_n \}_{ n \in \mathbb{ N } }$.
A draw from the prior is then generated by sampling $m \sim p_m$ and $n \sim q_n$, followed by $c | m, n \sim U( \Theta_c^{ ( m, n ) } )$.

Samples from this prior satisfy the four conditions listed at the top of this paragraph by construction, and it is also clear that condition \ref{cor_new} holds with $\Pi_c$-probability 1, leaving only condition \ref{cor_9} of Corollary \ref{explicit_conds} to verify.
To that end, fix $c_0 \in \Theta_c'$, a compact set $A \subset \Omega$ and $\varepsilon > 0$.
Fix $n' \in \N$ be such that $\varepsilon_{ n' } < \varepsilon$, and $m' \in \N$ such that $A \subseteq B_0( m' )$.
Then
\begin{equation*}
\Pi_c\left( c : \sup_{ \x \in A, \z \in J }\{ \| c( \x, \z ) - c_0( \x, \z ) \|_2 \} < \varepsilon \right) > \frac{ p_{ m' } q_{ n' } }{ | \Theta_c^{ ( m, n ) } | } > 0,
\end{equation*}
as required.

\subsection{Dirichlet process mixture model priors}\label{dirichlet_prior}

Let $\phi_{ \tau }( \mathbf{ z } )$ denote the $d$-dimensional centred Gaussian density with covariance matrix $\tau^{ -1 } \mathds{ I }_{ d \times d }$ restricted to $J$, as defined in Section \ref{net_prior}, and renormalised to be a probability density.
Let $F$ be a probability measure on $( 0, \infty )$ assigning positive mass to all non-empty open sets, and let $\operatorname{DP}( \zeta )$ denote the law of a Dirichlet process \citep{Ferguson73} with the mean measure $\zeta \in \M_f( J )$, which is taken to be a probability measure with a finite first moment, independent of $F$.
Let $\mathcal{ D }_{ \Upsilon }( J )$ denote the space of continuous, positive densities on $J$ with total mass at most $\Upsilon > 0$.
The Dirichlet process mixture model on $\mathcal{ D }_{ \Upsilon }( J )$ with truncated Gaussian mixture kernel $\phi_{ \tau }$ and mixing distribution $U( 0, \Upsilon ) \otimes F \otimes \operatorname{DP}( \zeta )$ is specified via the following sampling procedure:
\begin{enumerate}
\item Sample $P \sim \operatorname{DP}( \zeta )$. Then $P$ is a discrete probability measure on $\mathbb{ R }^d$ with countably many atoms with $\operatorname{DP}( \zeta )$-probability 1 \citep{Ferguson73}. Let $\z_1, \z_2, \ldots$ denote these atoms in some fixed ordering.
\item Sample IID copies $\tau_1, \tau_2, \ldots \sim F$.
\item Sample $\alpha \sim U( 0, \Upsilon )$.
\item Set $M( \z ) = \alpha \sum_{ j = 1 }^{ \infty } P( \z_j ) \phi_{ \tau_j }( \z - \z_j )$.
\end{enumerate}
Note that samples are finite measures with strictly positive, bounded densities on $J$.
Because $J$ is compact, $\mathcal{ D }_{ \Upsilon }( J )$ consists of continuous densities, and $F( ( x, \infty ) ) > 0$ for every $x > 0$, we also have that \citep[Theorem 1]{Bhattacharya12} holds, and hence the support of this prior is dense in $\mathcal{ D }_{ \Upsilon }( J )$.
Thus, conditions \ref{cor_11}, \ref{cor_13}, and \ref{cor_14} of Corollary \ref{explicit_conds} hold.

\section{Discussion}\label{discussion}

In this paper we have shown that posterior consistency for identifiable, nonparametric Bayesian inference of drift and jump coefficients of jump diffusions from discrete data holds under criteria which can be readily checked in practice.
Products of discrete net, wavelet, and Dirichlet process mixture model priors were shown to satisfy the conditions for consistency.

Our results shares the limitation of \citep{Gugushvili14}, in that we require priors to assign full mass to sets of functions for which the Lipschitz condition \eqref{lip} holds uniformly.
This rules out many widely used families of priors, such as Gaussian measures, but counterexamples exist to show that without it, uniform equicontinuity fails even for one dimensional unit diffusions [Matthias Birkner, personal communication].
It seems clear that an entirely different approach is needed, if consistency results are to be established without a uniform equicontinuity condition.

A further limitation of \citep{vanderMeulen13, Gugushvili14} is that of being established for a weak topology, for which the martingale approach of \citep{Walker04, Lijoi04} is well suited.
A testing approach, such as that of \citep{Ghosal07}, would yield convergence in a stronger topology as well as rates of convergence, but it is not clear how to adapt their results to the diffusion or jump diffusion settings.
Currently, results in this direction are only available for continuously observed scalar diffusions \citep{vanderMeulen06, Panzar09, Pokern13}, as well as discretely observed scalar diffusions on a compact interval \citep{Nickl17}.
However, these rely, respectively, on the continuity of the observation and on a tractable representation of the stationary density, neither of which is available in our setting.

Practical implementation of inference algorithms is beyond the scope of this paper, but we note that algorithms based on exact simulation for jump diffusions are available, at least in the scalar case \citep{Casella11, Goncalves11}.
Exact simulation of jump diffusions is an active area of research \citep{Goncalves13, Pollock15, Pollock16} and well suited for applications in Monte Carlo inference algorithms, with preliminary results in the continuous diffusion setting indicating that nonparametric algorithms can be feasibly implemented \citep{Papaspiliopoulos12, vanZanten13, vanderMeulen14}.
As a final remark, we note that presently such algorithms are only available for processes with jumps driven by compound Poisson processes of finite intensity, and with coefficients satisfying regularity assumptions comparable to those in Proposition \ref{existence}.
Thus our Theorem \ref{ts_consistency} brings the theory on nonparametric posterior consistency in line with current state of the art algorithms in one dimension, and anticipates development of comparable methods in higher dimensions.

\section*{Acknowledgements}

The authors are grateful to Matthias Birkner, the associate editor, and two anonymous referees for comments which have substantially improved this manuscript.
Jere Koskela was supported by EPSRC as part of the MASDOC DTC at the University of Warwick, grant No.~EP/HO23364/1, and by the Deutsche Forschungsgemeinschaft (DFG) grant BL 1105/3-2. 
Paul Jenkins is supported in part by EPSRC grant EP/L018497/1.

\bibliography{science}  
\bibliographystyle{plainnat}

\end{document}